\documentclass[10pt]{amsart} 
\usepackage{amscd}
\usepackage{amsmath,empheq}
\usepackage{amsfonts}
\usepackage{amssymb}
\usepackage{mathrsfs}
\usepackage[all]{xy}
\usepackage{mathtools}
\newtheorem{theorem}{Theorem}[section]

\newenvironment{proof-sketch}{\noindent{\bf Sketch of Proof}\hspace*{1em}}{\qed\bigskip}

\everymath{\displaystyle}

\renewcommand{\leq}{\leqslant}

\renewcommand{\geq}{\geqslant}
\newcommand{\ep}{\varepsilon}

\newcommand{\RR}{\mathbb R}

\newcommand{\ri}{\rightarrow}
\newcommand{\di}{\displaystyle}
\newcommand{\bb}{\begin{equation}}
\newcommand{\bbb}{\end{equation}}

\newcommand{\intom}{\int_\Omega}

\newcommand{\dpx}{\Delta_{p(x)}}

\newcommand{\calE}{{\mathcal E}}

\newcommand{\calX}{{\mathcal X}}
\newcommand{\wpp}{{W_0^{1,p(x)}(\Omega)}}

\begin{document}
\title[Non-autonomous singular biharmonic problems]
{Combined effects for non-autonomous \\ singular biharmonic problems}
\author[V.D. R\u{a}dulescu]{Vicen\c{t}iu D. R\u{a}dulescu}
\address[V.D. R\u{a}dulescu]{Faculty of Applied Mathematics, AGH University of Science and Technology, 30-059 Krak\'ow, Poland \& Institute of Mathematics, Physics and Mechanics, 1000 Ljubljana, Slovenia \& Department of Mathematics, University of Craiova, 200585 Craiova, Romania}
\email{\tt vicentiu.radulescu@imfm.si}
\author[D.D. Repov\v{s}]{Du\v{s}an D. Repov\v{s}}
\address[D.D. Repov\v{s}]{Faculty of Education and Faculty of Mathematics and Physics, University of Ljubljana \& Institute of Mathematics, Physics and Mechanics, 1000 Ljubljana, Slovenia}
\email{\tt dusan.repovs@guest.arnes.si}
\begin{abstract}
 We study the existence of nontrivial weak solutions for a class of generalized $p(x)$-biharmonic equations with singular nonlinearity and Navier boundary condition. The proofs combine variational and topological arguments. The approach developed in this paper allows for the treatment of several classes of singular biharmonic problems with variable growth arising in applied sciences, including the capillarity equation and the mean curvature problem.
\end{abstract}
\maketitle  
\thanks{\centerline{\it \small Dedicated to Professor Patrizia Pucci on the occasion of her 65th birthday} } 

\section{Introduction}
One of the numerous contributions of Patrizia Pucci to the development of the nonlinear analysis and mathematical physics concerns the refined qualitative analysis of real world phenomena driven by nonhomogeneous differential or nonlocal operators with one or more variable exponents. We refer only  to P.~Pucci and Q.~Zhang \cite{puccizh} for the existence of entire solutions for several classes of nonlinear problems with variable growth, G.~Autuori and P.~Pucci \cite{aupucci} for the study of the asymptotic stability for Kirchhoff systems in variable exponent Sobolev spaces, G.~Autuori, F.~Colasuonno and P.~Pucci \cite{aupucci1} for the analysis of stationary  higher-order Kirchhoff problems, and J.~Liu, P.~Pucci, H.~Wu, and Q.~Zhang \cite{liu} for the existence and blow-up rate of large solutions of $p(x)$-Laplacian equations with gradient terms.

The interest in the mathematical analysis of
partial differential equations driven by nonhomogeneous differential operators is motivated by their relevant applications in various disciplines. For instance, several models in the applied sciences are characterized by the fact that the associated energy density changes its ellipticity and growth properties according to the point. Such phenomena have been studied starting with the seminal works of T.C.~Halsey \cite{halsey} and V.V.~Zhikov \cite{zh1, zh2}, in close relationship with the qualitative mathematical analysis of strongly anisotropic materials
in the context of the homogenization and nonlinear elasticity. In the framework of materials with non-homogeneities, the standard approach based on the classical theory of $L^p$ and $W^{1,p}$ Lebesgue and Sobolev spaces is inadequate. We refer to electrorheological (smart) fluids or to phenomena in image processing, which should enable that the exponent $p$ is varying; see Y.~Chen, S.~Levine and M.~Rao \cite{clr}, and M.~Ruzicka \cite{ruzi}.
For instance, we refer to the Winslow effect of some fluids (like
lithium polymetachrylate) in which the viscosity in an electrical field is inversely proportional
to the strength of the field. The field induces string-like formations in the fluid,
which are parallel to the field. They can raise the viscosity by as much as five
orders of magnitude. This corresponds to electrorheological (non-Newtonian) fluids, which are mathematically described by means of nonlinear equations with variable exponent.
Such a study corresponds to the abstract setting of variable exponents Lebesgue and Sobolev spaces, $L^{p(x)}$ and $W^{1,p(x)}$, where $p$ is a real-valued function.

The theory of function spaces with variable exponent has been rigorously developed in the monograph of L.~Diening, P.~H\"asto,  P.~Harjulehto and M.~Ruzicka
 \cite{die} while the recent book by V.D.~R\u adulescu and D.D.~Repov\v{s} \cite{radrep} is devoted to the thorough variational and topological analysis of several classes of problems with one or more variable exponents; see also  the survey papers of P.~Harjulehto, P.~H\"ast\"o,  U.V.~Le and M.~Nuortio \cite{har1} and V.D.~R\u adulescu \cite{radnla}. We also refer to  G.~Mingione {\it et al.} \cite{mingi1, mingi2, mingi3},
 M.~Cencelj, V.D.~R\u adulescu and D.D.~Repov\v{s} \cite{cencel0}
 M.~Cencelj, D.D.~Repov\v{s} and \v{Z}.~Virk \cite{cencel}, and D.D.~Repov\v{s} \cite{repovs} for related results. The abstract setting of $p(x)$-biharmonic problems with singular weights has been recently considered by K.~Kefi and V.D.~R\u adulescu \cite{kefi} in relationship with microelectromechanical phenomena, surface diffusion on solids, thin film theory, flow in Hele-Shaw cells and phase field models
of multiphasic systems. The present paper extends  and complements  some results contained in \cite{kefi,kefi0}
to more general operators.

The study of  elliptic problems with variable exponent has recently been extended by I.H.~Kim and Y.H.~Kim \cite{kim} to a new class of non-homogeneous differential operators. Their contribution is a step forward in the analysis of nonlinear problems with variable exponent since it enables the understanding of problems with possible lack of uniform convexity. More precisely, in \cite{kim} they studied problems of the type
\bb\label{pro1}
\left\{\begin{array}{lll}
&\di -{\rm div}\, (\phi (x,|\nabla u|)\nabla u)=f(x,u)&\quad\mbox{in}\ \Omega\\
&\di u=0&\quad\mbox{on}\ \partial\Omega,
\end{array}\right.
\bbb
where $\Omega\subset\RR^N$ is a bounded domain with smooth boundary. The nonlinear term $f:\Omega\times\RR\ri\RR$ satisfies the Carath\'eodory condition and the function $\phi(x,t)$ is of type $|t|^{p(x)-2}$ with $p:\overline\Omega\ri (1,\infty)$ continuous.

 In the special case when $\phi(x,t)=|t|^{p(x)-2}$, the operator involved in problem \eqref{pro1} reduces to the $p(x)$-Laplacian, that is, $$\dpx u={\rm div}\, (|\nabla u|^{p(x)-2}\nabla u).$$ In many papers (see, e.g., \cite[Hypothesis (A4), p. 2629]{mrroyal}), the functional $\Phi$ induced by the principal part of problem \eqref{pro1} is assumed to be uniformly convex, namely, there exists $k>0$ such that for all $x\in\Omega$ and all $\xi,\psi\in\RR^N$,
$$\Phi\left(x,\frac{\xi+\psi}{2} \right)\leq\frac 12\, \Phi (x,\xi)+\frac 12\,\Phi (x,\psi)-k\,|\xi-\psi|^{p(x)}.$$
However, since the function $\Psi (x,t)=t^p$ is not uniformly convex for $t>0$ and $1<p<2$, this condition is not applicable to all $p$-Laplacian problems. An important feature of the abstract setting developed in \cite{kim} is that the main results are obtained without any uniform convexity assumption.

In the present paper we extend \cite{kim} to problems involving $p(x)$-biharmonic operators and we describe some qualitative properties in the presence of singular terms. We develop the study of biharmonic problems with Navier boundary condition for equations driven by the operator
$\Delta (\phi (x,|\Delta u|)\Delta u),$
where $\phi$ is as in \eqref{pro1}. Notice that if $\phi(x,t)=|t|^{p(x)-2}$, then we obtain the $p(x)$-biharmonic operator defined by $\Delta^2_{p(x)}u=\Delta(|\Delta u|^{p(x)-2}\Delta u).$

\section{Abstract Framework and Preliminary Results}\label{sec2}
Throughout this paper we assume that $\Omega\subset\RR^N$ is a bounded domain with smooth boundary.

Set
$$C_+(\overline\Omega)=\{h\in C(\overline\Omega),\;h(x)>1\;{\rm
for}\;
{\rm all}\;x\in\overline\Omega\}.$$
Assume that $p\in C_+(\overline\Omega)$ and let
$$p^+=\sup_{x\in\Omega}p(x)\qquad\mbox{and}\qquad p^-=
\inf_{x\in\Omega}p(x).$$

We define the Lebesgue space with variable exponent
 by
$$L^{p(x)}(\Omega)=\left\{u;\ u\ \mbox{is
 measurable and }
\int_\Omega|u(x)|^{p(x)}\;dx<\infty\right\}.$$
This function space is a Banach space if it is endowed with the norm
$$|u|_{p(x)}=\inf\left\{\mu>0;\;\int_\Omega\left|
\frac{u(x)}{\mu}\right|^{p(x)}\;dx\leq 1\right\}.$$
This norm is also called the Luxemburg norm.
Then $L^{p(x)}(\Omega)$ is reflexive if and only if $1 < p^-\leq p^+<\infty$
 and continuous functions with compact support
are dense in $L^{p(x)}(\Omega)$ if $p^+<\infty$.

The standard inclusion between Lebesgue spaces generalizes to the framework of spaces with variable exponent, namely if
$0 <|\Omega|<\infty$ and $p_1$, $p_2$ are variable exponents such that $p_1\leq p_2$
 in $\Omega$ then there exists the continuous embedding
$L^{p_2(x)}(\Omega)\hookrightarrow L^{p_1(x)}(\Omega)$.

Let $L^{p'(x)}(\Omega)$ denote the conjugate space
of $L^{p(x)}(\Omega)$, where $1/p(x)+1/p'(x)=1$. Then for all
$u\in L^{p(x)}(\Omega)$ and $v\in L^{p'(x)}(\Omega)$  the following H\"older-type inequality holds:
\begin{equation}\label{Hol}
\left|\int_\Omega uv\;dx\right|\leq\left(\frac{1}{p^-}+
\frac{1}{p'^-}\right)|u|_{p(x)}|v|_{p'(x)}\,.
\end{equation}

An important role in analytic arguments on Lebesgue spaces with variable exponent is played by the {\it modular} of  $L^{p(x)}(\Omega)$, which
is the map
 $\rho_{p(x)}:L^{p(x)}(\Omega)\rightarrow\RR$ defined by
$$\rho_{p(x)}(u)=\int_\Omega|u|^{p(x)}\;dx.$$
If $(u_n)$, $u\in L^{p(x)}(\Omega)$ and $p^+<\infty$ then the following
properties
hold:
\begin{equation}\label{L4}
|u|_{p(x)}>1\;\;\;\Rightarrow\;\;\;|u|_{p(x)}^{p^-}\leq\rho_{p(x)}(u)
\leq|u|_{p(x)}^{p^+}
\end{equation}
\begin{equation}\label{L5}
|u|_{p(x)}<1\;\;\;\Rightarrow\;\;\;|u|_{p(x)}^{p^+}\leq
\rho_{p(x)}(u)\leq|u|_{p(x)}^{p^-}
\end{equation}
\begin{equation}\label{L6}
|u_n-u|_{p(x)}\rightarrow 0\;\;\;\Leftrightarrow\;\;\;\rho_{p(x)}
(u_n-u)\rightarrow 0.
\end{equation}

We define the variable exponent  Sobolev space by
$$
W^{1,p(x)}(\Omega)=\{u\in L^{p(x)}(\Omega):\;|\nabla u|\in
L^{p(x)} (\Omega) \}.
$$

On $W^{1,p(x)}(\Omega)$ we may consider one of the following
equivalent norms
$$
\|u\|_{p(x)}=|u|_{p(x)}+|\nabla u|_{p(x)}
$$
or
$$\|u\|_{p(x)}=\inf\left\{\mu>0;\;\int_\Omega\left(\left| \frac{\nabla
u(x)}{\mu}\right|^{p(x)}+\left|
\frac{u(x)}{\mu}\right|^{p(x)}\right)\;dx\leq 1\right\}\,.
$$

Zhikov \cite{zh2} showed that smooth functions are in general not dense in
$W^{1,p(x)}(\Omega)$. This property is in relationship with the {\it Lavrentiev phenomenon}, which
asserts that there exist variational problems for which the infimum over the
smooth functions is strictly greater than the infimum over all functions that
satisfy the same boundary conditions. We refer to \cite[pp. 12-13]{radrep} for more details.

Let $W_0^{1,p(x)}(\Omega)$ denote the closure of the set of compactly supported $W^{1,p(x)}$-functions with respect to the norm $\|u\|_{p(x)}$.
When smooth functions are dense, we can also use the closure of
$C_0^\infty(\Omega)$ in $W^{1,p(x)}(\Omega)$. Using the Poincar\'e inequality, the space $W_0^{1,p(x)}(\Omega)$ can be  defined, in an equivalent manner, as the closure of $C_0^\infty(\Omega)$ with respect to the norm
$$\|u\|_{p(x)}=|\nabla u|_{p(x)}.$$

The vector space $(W^{1,p(x)}_0(\Omega),\|\cdot\|)$ is a separable and
reflexive Banach space.
Moreover, if $0 <|\Omega|<\infty$ and $p_1$, $p_2$ are variable exponents such that $p_1\leq p_2$
 in $\Omega$ then there exists a continuous embedding
$W^{1,p_2(x)}_0(\Omega)\hookrightarrow W^{1,p_1(x)}_0(\Omega)$.

Set
\bb\label{rho2}
\varrho_{p(x)}(u)=\int_\Omega |\nabla
u(x)|^{p(x)}\,dx.
\bbb

If $(u_n)$, $u\in W^{1,p(x)}_0(\Omega)$ then the following
properties hold:
\begin{equation}\label{M4}
\|u\|>1\;\Rightarrow\;\|u\|^{p^-}\leq \varrho_{p(x)}(u)
\leq\|u\|^{p^+}\,,
\end{equation}
\begin{equation}\label{M5}
\|u\|<1\;\Rightarrow\;\|u\|^{p^+}\leq \varrho_{p(x)}(u)
\leq\|u\|^{p^-}\,,
\end{equation}
\begin{equation}\label{M6}
\|u_n-u\|\rightarrow 0\;\;\;\Leftrightarrow\;\;\;\varrho_{p(x)}
(u_n-u)\rightarrow 0\,.
\end{equation}

Set
$$p_*(x)=\left\{
\begin{array}{lll}
&\di\frac{Np(x)}{N-p(x)}&\quad
\mbox{if $p(x)<N$}\\
&\di  +\infty&\quad\mbox{if $p(x)\geq N$}.
\end{array}\right.
$$
We point out that if $p,q\in C_+(\overline\Omega)$
and $q(x)<p_\star(x)$ for all $x\in\overline\Omega$ then the
embedding
$W_0^{1,p(x)}(\Omega)\hookrightarrow L^{q(x)}(\Omega)$
is compact.

For any positive integer $k$,  let $$W^{k,p(x)}(\Omega)=\{u\in L^{p(x)}(\Omega): D^{\alpha}u\in L^{p(x)}(\Omega), |\alpha|\leq k\},$$ where $\alpha=(\alpha_1,\alpha_2,...,\alpha_N) $ is a multi-index, $|\alpha|=\sum_{i=1}^{N}\alpha_i$ and
$$D^{\alpha}u=\frac{\partial^{|\alpha|}u}{\partial^{\alpha_1}x_1\ldots\partial^{\alpha_N}x_N}.$$
Then $W^{k,p(x)}(\Omega)$ is a separable and reflexive Banach space equipped with the norm $$\|u\|_{k,p(x)}=\displaystyle\sum_{|\alpha|\leq k}|D^{\alpha}u|_{p(x)}.$$

The space $W_{0}^{k,p(x)}(\Omega)$ is the closure of $C_{0}^{\infty}(\Omega)$ in $W^{k,p(x)}(\Omega).$

Next, we recall some properties of the space $$\calX := W^{1,p(x)}_{0}(\Omega)\cap W^{2,p(x)}(\Omega).$$
For any $u\in \calX$ we have $\|u\| = \|u\|_{1,p(x)}+\|u\|_{2,p(x)}$, thus  $$\|u\| =|u|_{p(x)}+|\nabla u|_{p(x)}+\displaystyle\sum_{|\alpha|=2}|D^{\alpha}u|_{p(x)}.$$

 A.B.~Zang and Y.~Fu \cite{ZF} proved the equivalence of the norms and they also established that the norm
$|\Delta u|_{p(x)}$ is equivalent to the norm $\|u\| $ (see \cite[Theorem 4.4]{ZF}).
 Note that $(\calX, \|.\|)$ is  a separable and
reflexive Banach space.

We recall that the critical Sobolev exponent is defined as follows:
 $$
\left\{
\begin{array}{l}
p^*(x)=\displaystyle\frac{Np(x)}{N-2p(x)},\quad p(x)<\frac{N}{2},\\
p^*(x)=\displaystyle+\infty, \hspace{0.8cm}\quad\quad p(x)\geq \frac{N}{2}. \\
\end{array}
\right.
$$

Assume that $q \in C^+(\overline{\Omega})$ and $q(x) < p^*(x)$ for any $x\in\Omega$. Then, by Theorem 3.2 in \cite{AA}, the function space $\calX$ is
continuously and compactly embedded in $L^{q(x)}(\Omega)$.

For a constant function $p$, the variable exponent Lebesgue and Sobolev
spaces coincide with the standard Lebesgue and Sobolev spaces. As pointed out in \cite{radrep}, the function spaces with variable exponent have some striking properties, such as:

(i)
If $1<p^-\leq p^+<\infty$ and $p:\overline\Omega\ri [1,\infty)$ is smooth, then the formula
$$\int_\Omega |u(x)|^pdx=p\int_0^\infty t^{p-1}\,|\{x\in\Omega ;\ |u(x)|>t\}|\,dt$$
has  no variable exponent analogue.

(ii) Variable exponent Lebesgue spaces do {\it not} have the {\it mean continuity property}. More precisely, if $p$ is continuous and nonconstant in an open ball $B$, then there exists a function $u\in L^{p(x)}(B)$ such that $u(x+h)\not\in L^{p(x)}(B)$ for all $h\in\RR^N$ with arbitrary small norm.

(iii) The function spaces with variable exponent
 are {\it never} translation invariant.  The use
of convolution is also limited, for instance  the Young inequality
$$| f*g|_{p(x)}\leq C\, | f|_{p(x)}\, \| g\|_{L^1}$$
holds if and only if
$p$ is constant.

\section{The Main Result}\label{sec3}
Throughout this paper we assume that $\Omega\subset\RR^N$ is a bounded domain with smooth boundary.

Let $p\in C_+(\overline\Omega)$ and consider the function $\phi:\Omega\times [0,\infty)\ri [0,\infty)$ satisfying the following hypotheses:

\smallskip
\noindent (H1) the mapping $\phi(\cdot,\xi)$ is measurable on $\Omega$ for all $\xi\geq 0$ and $\phi(x,\cdot)$ is locally absolutely continuous on $[0,\infty)$ for almost all $x\in\Omega$;

\smallskip\noindent (H2) there exist $a\in L^{p'}(\Omega)$ and $b>0$ such that
$$|\phi (x,|v|)v|\leq a(x)+b|v|^{p(x)-1}$$
for almost all $x\in\Omega$ and for all $v\in\RR^N$;

\smallskip\noindent (H3) there exists $c>0$ such that
$$\phi(x,\xi)\geq c\xi^{p(x)-2},\quad \phi(x,\xi)+\xi\frac{\partial\phi}{\partial\xi}(x,\xi)\geq c\xi^{p(x)-2}$$
for almost all $x\in\Omega$ and for all $\xi>0$.

\smallskip
An interesting consequence of theses assumptions is that $\phi$ satisfies a Simon-type inequality. More precisely, if we denote
$$\Omega_1:=\{x\in\Omega:\ 1<p(x)<2\}\quad\mbox{and}\quad\Omega_2:=\{x\in\Omega;\ p(x)\geq 2\},$$
then the following estimate holds for all $u,v\in\RR^N$
\bb\label{simonvari}\begin{array}{ll}
&\di \langle\phi (x,|u|)u-\phi (x,|v|)v,u-v\rangle\\
&\di\geq\left\{\begin{array}{lll}
&\di c(|u|+|v|)^{p(x)-2}|u-v|^2&\quad\mbox{if}\ x\in\Omega_1\ \mbox{and}\ (u,v)\not=(0,0)\\
&\di 4^{1-p^+}c|u-v|^{p(x)}&\quad\mbox{if}\ x\in\Omega_2,
\end{array}\right.
\end{array}
\bbb
 where $c$ is the positive constant from hypothesis (H3).

 Let $A:W_0^{1,p(x)}(\Omega)\ri\RR$ defined by
 $$A(u)=\intom\int_0^{|\nabla u(x)|}s\phi(x,s)dsdx.$$
  Inequality \eqref{simonvari}
 was used in \cite{kim} to show that $A':\wpp\ri W^{-1,p'(x)}(\Omega)$ is both a monotone operator and a mapping of type $(S_+)$.
We refer to Simon \cite{simon} for the initial version of inequality \eqref{simonvari} in the framework of the $p$-Laplace operator.

\smallskip
We study the following singular biharmonic problem with variable growth:
\bb\label{problem}
\left\{\begin{array}{lll}
&\di \Delta (\phi(x,|\Delta u|)\Delta u)=|u|^{-q(x)-1}u+\lambda |u|^{r(x)-2}u, \ & x\in \Omega\\
&\di u=\Delta u=0, \ & x\in \partial\Omega,
\end{array}\right.
\bbb
where $q$, $r$ are continuous functions and $\lambda$ is a positive parameter.

 If $\phi (x,\xi)=\xi^{p(x)-2}$ then we obtain the standard $p(x)$-Laplace biharmonic operator, that is, $\Delta^2_{p(x)}u:=\Delta (|\Delta u|^{p(x)-2}\Delta u)$.

Our abstract setting includes the case $\phi (x,\xi)=(1+|\xi|^2)^{(p(x)-2)/2}$, which corresponds to the generalized biharmonic mean curvature operator
$$\Delta \left[(1+|\Delta u|^2)^{(p(x)-2)/2}\Delta u \right].$$
The biharmonic capillarity equation corresponds to
$$\phi(x,\xi)=\left(1+\frac{\xi^{p(x)}}{\sqrt{1+\xi^{2p(x)}}}\right)\xi^{p(x)-2},\quad x\in\Omega,\ \xi>0,$$ hence the corresponding capillary phenomenon is described by the differential operator
$$\Delta \left[\left( 1+\frac{|\Delta u|^{p(x)}}{\sqrt{1+|\Delta u|^{2p(x)}}} \right)|\Delta u|^{p(x)-2}\Delta u\right].$$

We say that $u$ is a solution of problem \eqref{problem} if $u\in \calX\setminus\{0\}$ with $\Delta u=0$ on $\partial\Omega$ and
$$\intom \phi(x,|\Delta u|)\Delta u\Delta vdx=\intom |u|^{-q(x)-1}uvdx+\lambda\intom |u|^{r(x)-2}uv,$$
for all $v\in \calX$.

\begin{theorem}\label{t1}
Assume that hypotheses (H1)--(H3) are fulfilled and that
\bb\label{ipoteza}0<q(x)<1<r(x)<p(x)<p^*(x)\quad\mbox{for all}\ x\in\Omega.\bbb
Then problem \eqref{problem} has a solution for all $\lambda>0$.
\end{theorem}

In the present paper, problem \eqref{problem} is studied for the {\it subcritical case}, namely under the basic hypothesis \eqref{ipoteza}, which is crucial for compactness arguments. We consider that a very interesting research direction is to study the same problem in the {\it almost critical} setting, hence under the following assumption: there exists $x_0\in\Omega$ such that
\bb\label{critical} r(x)<p^*(x)\ \mbox{for all $x\in\Omega\setminus\{x_0\}$}\ \mbox{and}\ r(x_0)=p^*(x_0).\bbb
Of course, this hypothesis is not possible if the functions $p$ and $r$ are {\it constant}.
We conjecture that the result stated in Theorem \ref{t1} remains true under assumption \eqref{critical}.

\section{Proof of Theorem \ref{t1}} Fix $\lambda>0$ and denote $$\Phi(x,t):=\int_0^ts\phi(x,s)ds\quad\mbox{for all}\ x\in\Omega.$$

The energy functional associated to problem \eqref{problem} is $\calE:\calX\ri\RR$ defined by
$$\calE (u)=\intom\Phi(x,|\Delta u|)dx-\intom\frac{|u|^{1-q(x)}}{1-q(x)}dx-\lambda\intom\frac{|u|^{r(x)}}{r(x)}dx.$$
By hypothesis \eqref{ipoteza}, we deduce that $\calE$ is well-defined. On the other hand, with the same arguments as in \cite[Proposition 3.3]{kefi},
the energy functional $\calE$ is sequentially lower semicontinuous and of class $C^1$. Moreover, the mapping $\calE':\calX\ri\calX^*$ is a strictly monotone, bounded homeomorphism and is of type $(S_+)$, that is, if
$$u_n\rightharpoonup u\ \mbox{in}\ \calX\ \mbox{and}\ \limsup_{n\ri\infty}\calE'(u_n)(u_n-u)\leq 0,$$
then $u_n\ri u$ in $\calX$.

We split the proof of Theorem \ref{t1} into several steps.

\smallskip{\it Step 1.} The functional $\calE$ is coercive.

Using (H3), we first deduce that for all $u\in\calX$
$$\calE(u)\geq c\intom\frac{|\Delta u|^{p(x)}}{p(x)}dx-\intom\frac{|u|^{1-q(x)}}{1-q(x)}dx-\lambda\intom\frac{|u|^{r(x)}}{r(x)}dx.$$
Therefore
$$\calE(u)\geq \frac{c}{p^+}\intom |\Delta u|^{p(x)}dx-\frac{1}{1-q^+}\intom |u|^{1-q(x)}dx-\frac{\lambda}{r^-}\intom |u|^{r(x)}dx.$$
 It follows that for all $u\in\calX$ with $\|u\|>1$ we have
$$\begin{array}{ll}
\calE(u)&\geq\di \frac{c}{p^+}\,\|u\|^{p^-}-\frac{1}{1-q^+}\intom |u|^{1-q(x)}dx-\frac{\lambda}{r^-}\intom |u|^{r(x)}dx\\
&\di\geq \frac{c}{p^+}\,\|u\|^{p^-}-\frac{1}{1-q^+}\min\left\{|u|^{1-q^+}_{(1-q(x))p(x)},|u|^{1-q^-}_{(1-q(x))p(x)} \right\}\\
&\di - \frac{\lambda}{r^-}\min\{|u|^{r^+}_{r(x)},|u|^{r^-}_{r(x)} \}.
\end{array}$$
Next, by hypothesis \eqref{ipoteza}, it follows that there exists $c_0>0$ such that for all $u\in\calX$
$$ \max\{|u|_{(1-q(x))p(x)},|u|_{r(x)}\}\leq c_0\,\|u\|.$$
We deduce that
$$\begin{array}{ll}
\calE(u)&\geq\di \frac{c}{p^+}\,\|u\|^{p^-}-\frac{c_0}{1-q^+}\min\left\{\|u\|^{1-q^+},\|u\|^{1-q^-} \right\}
 - \frac{\lambda c_0}{r^-}\min\{\|u\|^{r^+},\|u\|^{r^-} \}\\
 &\di = \frac{c}{p^+}\,\|u\|^{p^-}-\frac{c_0}{1-q^+}\,\|u\|^{1-q^+}
 - \frac{\lambda c_0}{r^-}\,\|u\|^{r^-} \,.
\end{array}$$

We conclude the proof of Step 1 by using hypothesis \eqref{ipoteza}.

\smallskip
The next step shows that the energy $\calE$ does not satisfy one of the geometric hypotheses of the mountain pass theorem. More precisely, we show that there exists a ``valley" for $\calE$ close to the origin, so not far away from the origin, as it is required by the Ambrosetti-Rabinowitz theorem.

\smallskip{\it Step 2.} There exists $v\in\calX$ such that $\calE(tv)<0$ for all small enough $t>0$.

Hypothesis (H2) yields for all $u\in\calX$
$$\Phi(x,|\Delta u|)\leq\left|\int_0^{|\Delta u|}\left( a(x)+b|s|^{p(x)-1}\right)ds\right|\leq |a(x)|\,|\Delta u|+b\,\frac{|\Delta u|^{p(x)}}{p(x)}\,.$$
It follows that
$$\intom \Phi(x,|\Delta u|)dx\leq\intom |a(x)|\,|\Delta u|dx+b\intom \frac{|\Delta u|^{p(x)}}{p(x)}dx.$$

Fix $v\in C^\infty(\Omega)$ with $\mbox{supp}\, (v)\subset\Omega$ and $0\leq v\leq 1$. For all $t>0$ we have
$$\begin{array}{ll}
\di\intom \Phi(x,|\Delta (tv)|)dx&\di\leq 2t|a|_{p'(x)}\,|\Delta v|_{p(x)}+b\intom t^{p(x)}\,\frac{|\Delta v|^{p(x)}}{p(x)}dx\\
&\di\leq 2t|a|_{p'(x)}\,|\Delta v|_{p(x)}+b\,\frac{t^{p^-}}{p^-}\intom |\Delta v|^{p(x)}dx\\
&\di =C_1t+C_2t^{p^-},\end{array}$$
where $$C_1=2|a|_{p'(x)}\,|\Delta v|_{p(x)}>0\quad \mbox{and}\quad C_2=\frac{b}{p^-}\intom |\Delta v|^{p(x)}dx>0.$$

We conclude that
\bb\label{newdeal}\calE(tv)\leq C_1t+C_2t^{p^-}-C_3t^{1-q^-},\bbb
where $$C_3=\intom\frac{|v|^{1-q(x)}}{1-q(x)}dx>0.$$
Since $0<1-q^-<1<p^-$, relation \eqref{newdeal} implies that $\calE(tv)<0$, provided that $t>0$ is small enough.

\smallskip{\it Step 3.} The infimum of $\calE$ is achieved by some $u_0\in\calX\setminus\{0\}$.

Let $(u_n)\subset\calX$ be a minimizing sequence of $\calE$. By Step 1, we deduce that $(u_n)$ is a bounded sequence. So, there exists $u_0\in\calX$ such that, up to a subsequence,
$$u_n\rightharpoonup u_0\quad\mbox{in}\ \calX$$
$$u_n\ri u_0\quad\mbox{in}\ L^{r(x)}(\Omega).$$

By the weak lower semicontinuity of $\calE$ we conclude that
$$m:=\inf\{\calE(u);\ u\in\calX\}\leq \calE(u_0)\leq\liminf_{n\ri\infty}\calE(u_n)=m,$$
hence $u_0$ is a global minimizer of $\calE$
and
$$m=\calE(u_0).$$
Moreover, by Step 2, we have $m<0$, hence $u_0\in\calX\setminus\{0\}$.

\smallskip
To complete the proof of Step 3, it remain to show that $u_0$ satisfies \eqref{problem} in the weak sense and that $\Delta u_0=0$ on $\partial\Omega$.
These properties will be established in the final steps of the proof.

\smallskip{\it Step 4.} We have
\bb\label{caracal}\intom \phi(x,|\Delta u_0|)\Delta u_0\Delta vdx=\intom |u_0|^{-q(x)-1}u_0vdx+\lambda\intom |u_0|^{r(x)-2}u_0v,\bbb
for all $v\in \calX$.

Fix $v\in\calX$ and $\ep>0$. Define $z=(u_0+\ep v)^+$.

Since $u_0$ is a (global) minimizer of $\calE$, we deduce that
$$\begin{array}{ll}
0&\di\leq \intom \phi(x,|\Delta u_0|)\Delta u_0\Delta zdx-\intom |u_0|^{-q(x)-1}u_0zdx-\lambda\intom |u_0|^{r(x)-2}u_0zdx\\
&\di= \int_{[u_0+\ep v>0]} \phi(x,|\Delta u_0|)\Delta u_0\Delta (u_0+\ep v)dx-\int_{[u_0+\ep v>0]} |u_0|^{-q(x)-1}u_0(u_0+\ep v)dx\\
&\di-\lambda\int_{[u_0+\ep v>0]}|u_0|^{r(x)-2}u_0(u_0+\ep v)\\
&\di =\intom \phi(x,|\Delta u_0|)\Delta u_0\Delta (u_0+\ep v)dx-\intom |u_0|^{-q(x)-1}u_0(u_0+\ep v)dx\\
&\di-\lambda\intom |u_0|^{r(x)-2}(u_0+\ep v)dx\\
&\di -\int_{[u_0+\ep v\leq 0]} \phi(x,|\Delta u_0|)\Delta u_0\Delta (u_0+\ep v)dx+\int_{[u_0+\ep v\leq 0]} |u_0|^{-q(x)-1}u_0(u_0+\ep v)dx\\
&\di+\lambda\int_{[u_0+\ep v\leq 0]}|u_0|^{r(x)-2}u_0(u_0+\ep v)dx.
\end{array}
$$
It follows that
$$\begin{array}{ll}
0&\di\leq \intom \phi(x,|\Delta u_0|)|\Delta u_0|^2dx-\intom |u_0|^{1-q(x)}dx-\lambda\intom |u_0|^{r(x)}dx\\
&\di +\ep \intom \phi(x,|\Delta u_0|)\Delta u_0\Delta  vdx-\ep\intom |u_0|^{-q(x)-1}u_0 vdx\\
&\di-\lambda\ep\intom |u_0|^{r(x)-2}u_0 vdx+O(\ep^2)\\
&\di -\int_{[u_0+\ep v\leq 0]} \phi(x,|\Delta u_0|)\Delta u_0\Delta (u_0+\ep v)dx+\int_{[u_0+\ep v\leq 0]} |u_0|^{-q(x)-1}u_0(u_0+\ep v)dx\\
&\di+\lambda\int_{[u_0+\ep v\leq 0]}|u_0|^{r(x)-2}u_0(u_0+\ep v)dx.
\end{array}
$$

Since $u_0$ is a critical point of $\calE$, this relation yields
$$\begin{array}{ll}
0&\di\leq\ep\left( \intom \phi(x,|\Delta u_0|)\Delta u_0\Delta  vdx-\intom |u_0|^{-q(x)-1}u_0 vdx
-\lambda\intom |u_0|^{r(x)-2}u_0 vdx\right)\\
&\di -\ep\int_{[u_0+\ep v\leq 0]} \phi(x,|\Delta u_0|)\Delta u_0\Delta  vdx\\
&\di=\ep\left( \intom \phi(x,|\Delta u_0|)\Delta u_0\Delta  vdx-\intom |u_0|^{-q(x)-1}u_0 vdx
-\lambda\intom |u_0|^{r(x)-2}u_0 vdx\right)\\
&\di +o(\ep)\quad\mbox{as}\ \ep\ri 0.
\end{array}$$
This relation implies that
$$ \intom \phi(x,|\Delta u_0|)\Delta u_0\Delta  vdx-\intom |u_0|^{-q(x)-1}u_0 vdx
-\lambda\intom |u_0|^{r(x)-2}u_0 vdx\geq 0.$$
Changing $v$ with $-v$ we deduce that
$$ \intom \phi(x,|\Delta u_0|)\Delta u_0\Delta  vdx-\intom |u_0|^{-q(x)-1}u_0 vdx
-\lambda\intom |u_0|^{r(x)-2}u_0 vdx\leq 0.$$
We conclude that relation \eqref{caracal} holds.

\smallskip{\it Step 5.} We have $\Delta u_0=0$ on $\partial\Omega$.

We use some ideas developed in \cite[pp. 79-80]{kefi}.
By virtue of \eqref{caracal}, the function $u_0$ satisfies for all $v\in\calX$
\bb\label{15bis}
\intom \phi(x,|\Delta u_0|)\Delta u_0\Delta vdx=\intom A(x)vdx,\bbb
where
$$A(x):=|u_0|^{-q(x)-1}u_0+\lambda\,|u_0|^{r(x)-2}u_0.$$

Let $z\in\calX$ be the unique solution of the linear problem
\bb\label{16bis}
\left\{
\begin{array}{lll}
&\di \Delta z=A(x)&\quad\mbox{in}\ \Omega\\
&\di z=0&\quad\mbox{on}\ \partial\Omega.
\end{array}\right.
\bbb
Relations \eqref{15bis} and \eqref{16bis} yield that for all $v\in\calX$
$$\intom \phi(x,|\Delta u_0|)\Delta u_0\Delta vdx=\intom (\Delta z)vdx.$$
By Green's formula we deduce that for all $v\in C^\infty_c(\Omega)\subset\calX$
\bb\label{17bis}
\intom \phi(x,|\Delta u_0|)\Delta u_0\Delta vdx=\intom z\Delta vdx.\bbb

For any $w\in C^\infty_c(\Omega)$, let $v\in C^\infty_c(\Omega$ be the unique solution of the problem
$$\left\{
\begin{array}{lll}
&\di \Delta v=w&\quad\mbox{in}\ \Omega\\
&\di v=0&\quad\mbox{on}\ \partial\Omega.
\end{array}\right.
$$
Returning to \eqref{17bis}, we deduce that for all $w\in C^\infty_c(\Omega)$
$$\intom\left(\phi(x,|\Delta u_0|)\Delta u_0-z \right)wdx=0.$$
Applying \cite[Lemma VIII.1]{brezis} we conclude that
\bb\label{mitica}\phi(x,|\Delta u_0|)\Delta u_0-z=0\quad\mbox{in}\ \Omega.\bbb
But $z=0$ on $\partial\Omega$. Using hypothesis (H3), relation \eqref{mitica} implies that $\Delta u_0=0$ on $\partial\Omega$. The proof of Theorem \ref{t1} is now complete. \qed

We point out that the same arguments are no longer valid if the parameter $\lambda$ in problem \eqref{problem} is negative. In this case, the conclusion of Step~2 is not true, hence it is possible that the global minimizer $u_0$ obtained in Step~3 is trivial. Thus, if $\lambda<0$, the reaction term $|u|^{r(x)-2}u$ becomes a source term. We believe that if $\lambda$ is negative, a natural assumption is to replace the nonlinearity $|u|^{r(x)-2}u$ with a term having a different growth near the origin and at infinity.

\section{Epilogue} A very interesting open problem concerns the same analysis if the left-hand side  of problem \eqref{problem}
is replaced either by the differential operator
\bb\label{ene1} \Delta (\phi_1(x,|\Delta u|)\Delta u)+V(x) \Delta (\phi_2(x,|\Delta u|)\Delta u)\bbb
or by
\bb\label{ene2} \Delta (\phi_1(x,|\Delta u|)\Delta u)+V(x) \Delta (\phi_2(x,|\Delta u|)\Delta u)\log (e+|x|),\bbb
where $V$ is a nonnegative potential and $\phi_1$, $\phi_2$ satisfy hypotheses (H1)--(H3) corresponding to the variable exponents $p_1(x)$, $p_2(x)$ with  $p_1(x)\leq p_2(x)$ in $\Omega$.
Considering two different materials with power hardening exponents $p_1(x)$ and $p_2(x)$, respectively, the coefficient
$V(x)$ dictates the geometry of a composite of the
two materials. When $V(x)>0$, then
$p_2(x)$-material is present, otherwise the $p_1(x)$-material is the only
one making the composite.
  Composite materials with locally different hardening exponents
$p_1(x)$ and $p_2(x)$
can be described using the energies associated to the differential operators defined in \eqref{ene1} and \eqref{ene2}. 

Problems of this type were also motivated
by applications to elasticity, homogenization, modelling of strongly anisotropic materials, Lavrentiev phenomenon, etc.
 In the case of constant exponents, we refer to the pioneering papers by P.~Marcellini \cite{marce1, marce2} and G.~Mingione {\it et al.} \cite{mingi1, mingi2, mingi3}. Double phase problems with variable growth have recently been considered by M.~Cencelj, V.D.~R\u adulescu and D.D.~Repov\v{s} \cite{cencel0}, V.D.~R\u adulescu and Q.~Zhang \cite{radz}, and X.~Shi, V.D.~R\u adulescu, D.D.~Repov\v{s} and Q.~Zhang \cite{shi}.

 We conclude by pointing out that the differential operator
$\Delta (\phi(x,|\Delta u|)\Delta u)$
considered in problem \eqref{problem} falls in the realm of those related to the so-called Musielak-Orlicz spaces (see J.~Musielak \cite{musi} and W.~Orlicz \cite{orli}), more generally, of the operators having non-standard growth
conditions (which are widely considered in the calculus of variations). These function spaces are Orlicz spaces
whose defining Young function exhibits an additional dependence on the
$x$
variable. Nonlinear problems in Musielak-Orlicz spaces were studied in V.D.~R\u adulescu and D.D.~Repov\v{s} \cite[Chapter 4]{radrep}, but only in the framework of the standard $p(x)$-Laplace operator.

\medskip
{\bf Acknowledgements.} This research was supported by the Slovenian Research Agency grants
P1-0292, J1-8131, J1-7025, N1-0064, and N1-0083. V.D.~R\u adulescu acknowledges
also  the support through a grant from the Romanian Ministry of Research and Innovation, CNCS--UEFISCDI, project number PN-III-P4-ID-PCE-2016-0130,
within PNCDI III.

\end{document}